\documentclass[leqno]{article}

\usepackage{amssymb}
\usepackage{amsfonts}
\usepackage{latexsym}
\usepackage[T1,T2A]{fontenc}
\usepackage[cp1251]{inputenc}
\usepackage[T2B]{fontenc}
\usepackage[bulgarian,russian,english]{babel}

\newtheorem {theor} {\bf Theorem}
\newtheorem {lemma} {\bf Lemma}

\newtheorem{crl} {\bf Corollary}

\newtheorem{dfn}{\bf Definition}
\newtheorem {algr} {\bf Algorithm}

\title {Mathematical Modeling of the Weaving Structure Design
 }
\author {Krasimir Yordzhev\and Hristina Kostadinova}
\date {}

\begin {document}
\inputencoding{cp1251}

\maketitle
\begin{abstract}
  An equivalence relation in the set of all square binary matrices is described in this work. It is discussed a combinatoric problem about finding the cardinal number and the elements of the factor set according to this relation. We examine the possibility to get some special elements of this factor set. We propose an algorithm, which solves these problems. The results we have received are used to describe the topology of the different weaving structures.
\end{abstract}

\section{Introduction.} The present work demonstrates the ideas in the works \cite{umb2009} and \cite{umb2010} and in that sense it is their natural development.

As we know \cite{borzunob,yordzhev}, the interweaving of the fibres in certain weaving structure can be coded using square \textit{binary} (or (0,1), or  \textit{boolean}) matrix, i.e. all elements of this matrix are 0 or 1. The fabric represented by this matrix exists if and only if in each row and in each column of the matrix there is at least one zero and at least one one. Two different matrices correspond to one and the same weaving structure if and only if one matrix is made from the other one by several consecutive cycle moves of the first row or column to the last place.

Let $n$ is a whole positive number. Let denote ${\cal B}_n$ the set of all $n\times n$ binary matrices, and ${\cal Q}_n$ the set of all $n\times n$ binary matrices which have at least one one and one zero in every row and every column. It is obvious, that ${\cal Q}_n \subset {\cal B}_n$. About the necessary definitions and marks in the theory of matrices we conform with \cite{kurosh} and \cite{tarakanov}. It is not difficult to see, that
\begin{equation} \label{bn}
\left| {\cal B}_n \right| =2^{n^2}
\end{equation}

If $A=(a_{i\, j} ) \in {\cal B}_n$, then $A^T =(a_{j\, i} )$, $1\le i,j \le n$ denotes the {transposed} matrix $A$ \textit .

We are interested in the subset ${\cal P}_n \subset {\cal Q}_n$ made of all \textit{permutating} matrices, i.e. binary matrices which have exactly one one in every row and column. As it is well known \cite{tarakanov} the set ${\cal P}_n$ together with the operation multiplication of matrices is a group, isomorphic to the symmetric group ${\cal S}_n$, where the set
\begin{equation}\label{S_n}
{\cal S}_n =\left\{ \left. \left(
\begin{array}{cccc}
1    &    2    &  \cdots &    n\\
i_1  &    i_2  &  \cdots &    i_n
\end{array}
\right) \; \right| \; 0\le i_k \le n,\; k=1,2,\ldots ,n,\; i_k \ne i_l \; \mbox{\textrm{for}} \; k\ne l
\right\}
\end{equation}
is made of all one to one maps of the elements of the set $\{ 1,2,\ldots ,n\}$ to itself. And if $M\in {\cal P}_n$ and the corresponding element in this isomorphism is
$\displaystyle \left(
\begin{array}{cccc}
1    &    2    &  \cdots &    n\\
i_1  &    i_2  &  \cdots &    i_n
\end{array}
\right)
\in {\cal S}_n$, this means that the only one one in the first row of  $M$ to be on the $i_1$ -st place, the one in the second row of $M$ to be on the $i_2$-nd place, and so on, the one in the $n$th row of $M$ to be on the $i_n$th place.

Let $t\in \{1,2,\ldots ,n\}$ and let
$\displaystyle \rho =
\left(
\begin{array}{cccccc}
1    &    2    &  \cdots &  t    & \cdots  &  n  \\
i_1  &    i_2  &  \cdots &  i_t  & \cdots  &  i_n
\end{array}
\right)
\in {\cal S}_n$. We denote $(t)\rho =i_t$ the image $i_t$ of the number $t$ using the map $\rho$. However, for random $\rho_1 ,\rho_2 \in S_n$ by definition $(t)\rho_1 \rho_2 =((t)\rho_1 )\rho_2$ (see \cite{kurosh}).

As it is well-known \cite{tarakanov}, if we multiply random $n\times n$ matrix $A$ from the right with random permutational matrix $M\in {\cal P}_n$, then this is the same as changing the columns of $A$. And if the corresponding element of  $M\in {\cal P}_n$ in the above-described isomorphism is
$\displaystyle \left(
\begin{array}{cccc}
1    &    2    &  \cdots &    n\\
i_1  &    i_2  &  \cdots &    i_n
\end{array}
\right)
\in {\cal S}_n$,
then after the multiplication we get a matrix with $k$ column equal to $i_k$ column of $A$, $k=1,2,\ldots ,n$. Analogously when we want to exchange the rows we multiply $A$ from the left with $M^T$.

Identity element of the group ${\cal P}_n $ is the identity matrix $E_n$, consisting of ones in the leading diagonal and zeros everywhere else. The identity element of the group ${\cal S}_n$ is the element $\displaystyle \left(
\begin{array}{cccc}
1    &    2    &  \cdots &    n\\
1    &    2    &  \cdots &    n
\end{array}
\right)$,

We say, that the binary $n\times n$ matrices $A$ and $B$ are equivalent and we write $A\sim B$, if one matrix is made by the other after several consecutive cycle moves of the first row or column to the last place. In other words, if
$A,B\in {\cal Q}_n$ and $A\sim B$, then with the help of these matrices we code one and the same weaving structure (fabric). It is obvious, that the relation in the set ${\cal B}_n$ is a equivalence relation. The equivalence class according to the relation $\sim$ with the matrix $A$ we denote $\overline{A} $, and the sets of equivalence classes in ${\cal B}_n$  and   ${\cal Q}_n$ (factor set) according to $\sim$ with $\overline{{\cal B}_n}$ and $\overline{{\cal Q}_n}$. We consider that $\overline{{\cal B}_n}$  and $\overline{{\cal Q}_n}$ are described if there is a random representative of each equivalence class.
The equivalence classes of ${\cal B}_n$ by the equivalence  relation
$\sim$ are particular kind of  {\it double coset} (see \cite{hall} \S1.7, or
\cite{melnikov} v. 1, ch. 2, \S1.1). They make use of substitution
groups theory (see \cite{huppert} \S1.12, \S2.6) and linear
representation of finite groups theory (see \cite{curtis}
\S\S44-45).

The elements of $\overline{{\cal Q}_n}$ we call \textit{interweavings}. In that case the number $n$ is called \textit{repeating} of the interweavings of $\overline{{\cal Q}_n}$. These terms are taken from the Interweaving-knowing - the science which examines the design, physical and mechanical properties of the different interweavings of the fibres after the given textile structure is weaved.

It is naturally to arise a lot of combinatoric problems, which take place in practice in the weaving industry, connected with the different subsets of $\overline{{\cal Q}_n}$, i.e. with the different classes of interweavings. Some of these classes we examine in the present work.

\section{Some classes of interweavings.}
It is easy to see, that if $A\in {\cal P}_n$ and $B\sim A$, then $B\in {\cal P}_n$. Interweavings which representatives are elements of the set ${\cal P}_n$ of all permutational matrices are called \textit{primary} interweavings. A formula and an algorithm to calculate the number of all primary interweavings with a random repetition $n$ are shown in \cite{bansko,yordzhev}.

We examine the matrix
\begin{equation}\label{P}
P =
\left(
\begin{array}{ccccc}
0 & 1 &  0   & \cdots & 0 \\
0 & 0 &  1   & \cdots & 0 \\
\vdots    & \vdots    & \vdots   & \vdots  & \vdots   \\
0 & 0 &  0 & \cdots & 1  \\
1 & 0 &  0 & \cdots & 0  \\
\end{array}
\right)
=\left( p_{i\, j} \right) \in {\cal P}_n ,
\end{equation}
where $p_{1\, 2} = p_{2\, 3} = \cdots =p_{i\, i+1} = \cdots =p_{n-1\, n} = p_{n\, 1} = 1$ and these are the only one ones in $P$, and all other elements are zeros. At the above described isomorphism of the group of permutational matrices with the symmetric group, the matrix $P$ corresponds to the element
\begin{equation}\label{pi}
 \pi=
\left(
\begin{array}{cccccc}
1 & 2 & 3 & \cdots & n-1  &  n\\
2 & 3 & 4 & \cdots & n    &  1
\end{array}
\right)
\in {\cal S}_n .
\end{equation}

It is not difficult to calculate that
$P^t \ne E_n $ when $1\le t<n$, $P^n =E_n $,
where $E_n$ is identity matrix of row $n\times n$ and
$P^{k+n} =P^k $
for every natural number $k$.

Let $A\in {\cal B}_n$ and let
$$B=PA$$
and
$$C=A P .$$

It is easy to prove \cite{tarakanov}, that the first row of $B$ is equal to the second row of $A$, the second row of $B$ is equal to the third row of $A$ and so on, the last row of $B$ is equal to the first row of $A$, i.e. the matrix $B$ is created by the matrix $A$ by moving the first row to the last place, and the other rows are moved one level upper.

Analogously we convince that $C$ is made from $A$ by moving the last column to the first place, and the other column move one position to the right.

Have in mind what is above-described, it is easy to prove the following:

\begin{lemma}\label{l1}
Let  $A,B\in {\cal B}_n$. Then $A\sim B$,  if and only if, there exist natural numbers $k,l$, so that
$$A=P^k B P^l ,$$
 where $P$ is the matrix given by the formula (\ref{P}). From the equality  $P^{n+t} =P^t$  follows that for each natural number $t$, the numbers $k$ and $l$  is enough to be searched for in the interval $[0,n-1]$.

\hfill $\Box$
\end{lemma}

\begin{crl}\label{c1}
Each equivalence class of ${\cal B}_n$ according to the relation $\sim$ contains no more than $n^2$ elements.

\hfill $\Box$
\end{crl}

\begin{crl}\label{c2}
All elements of the given equivalence class of ${\cal B}_n$ according to the relation $\sim$ can be placed in a rectangular table with dimensions  $s\times t$, where $s$ and $t$ are divisors of $n$.

\hfill $\Box$
\end{crl}

Except the well-known from the linear algebra operation matrix transpose, here we take more similar matrix operations.

Let
\begin{equation}\label{A}
A=
\left(
\begin{array}{cccc}
a_{1\, 1} & a_{1\, 2} & \cdots & a_{1\, n} \\
a_{2\, 1} & a_{2\, 2} & \cdots & a_{2\, n} \\
\vdots    & \vdots    &        & \vdots    \\
a_{n\, 1} & a_{n\, 2} & \cdots & a_{n\, n} \\
\end{array}
\right)
\in {\cal B}_n
\end{equation}

If $A$ is a square binary matrix, represented by the formula (\ref{A}), then by definition
\begin{equation}\label{A^S}
A^S =
\left(
\begin{array}{cccc}
a_{1\, n} & a_{1\, n-1} & \cdots & a_{1\, 1} \\
a_{2\, n} & a_{2\, n-1} & \cdots & a_{2\, 1} \\
\vdots    & \vdots    &        & \vdots    \\
a_{n\, n} & a_{n\, n-1} & \cdots & a_{n\, 1} \\
\end{array}
\right)
,
\end{equation}
i. e. $A^S$  is created from $A$ as the last column of $A$ becomes the first, the column before last - second and so on, the first column becomes last. In other words, if $A=(a_{i\, j} )$, then $A^S =(a_{i\, n-j+1} )$, $1\le i,j \le n$.

It is obvious, that
$$ \left( A^S \right)^S =A.$$

We say, that the matrix $A\in{\cal B}_n$ is \textit{ a mirror image} of the matrix $B\in {\cal B}_n$, if $A^S =B$.

It is easy to see, that if the matrix $A$ is a mirror image of the matrix $B$, then $B$ is a mirror image of $A$, i.e. the relation ''mirror image'' is symmetric.

In the common case $A\ne A^S$. If $A = B^S$ and $B = C^S$, then in the common case we have
$A=B^S =\left( C^S \right)^S =C\ne C^S$.
Therefore, the relation ''mirror image'' is not reflexive and is not transitive.

We examine the matrix
\begin{equation}\label{S}
S =
\left(
\begin{array}{cccccc}
0 & 0 & \cdots & 0 & 0 & 1 \\
0 & 0 & \cdots & 0 & 1 & 0 \\
0 & 0 & \cdots & 1 & 0 & 0 \\
\vdots    & \vdots    &        & \vdots   & \vdots  & \vdots   \\
1 & 0 & \cdots & 0 & 0 & 0 \\
\end{array}
\right)
=\left( s_{i\, j} \right) \in {\cal P}_n ,
\end{equation}
where for each $i=1,2,\ldots ,n$ $s_{i\, n-i+1} =1$ и $s_{i\, j} =0$ as $j\ne n-i+1$. According to the above-described isomorphism of the group of permutational matrices with symmetric group, the matrix $S$ corresponds with the element
\begin{equation}\label{sigma}
\sigma =
\left(
\begin{array}{ccccc}
1 & 2 & 3 & \cdots & n\\
n & n-1 & n-2 & \cdots & 1
\end{array}
\right)
\in {\cal S}_n .
\end{equation}

Obviously $S$ is a symmetric matrix, i.e. $S^T = S$. We check directly, that $S^2 =E_n $.

It is not tough to notice (see for example \cite{tarakanov}), that for each $A\in {\cal B}_n$ is realized
\begin{equation}\label{8}
A^S = A S .
\end{equation}

\begin{lemma}\label{l2}
If $P$ and $S$ are the matrices given by the formulas (\ref{P}) and (\ref{S}), then for each $l=0,1,2,\ldots ,n-1$ the following proposition is true:
\begin{equation}\label{PlS}
P^l S=SP^{n-l}
\end{equation}
\end{lemma}

Proof.
Let denote $\oplus$ and $\ominus$ the operations corresponding to addition and subtraction in the ring ${\cal Z}_n =\{1,2,\ldots ,n\equiv 0 \}$ of the remainders by modulus $n$.
In order to  have accordance with the marks of the elements of ${\cal S}_n$ and as $0\equiv n\; (\mbox{\rm mod}\; n)$, then the zero element in ${\cal Z}_n$ we denote $n$ (instead of 0). If $\pi \in {\cal S}_n$ and $\sigma\in {\cal S}_n$ are elements corresponding to the matrices $P\in {\cal P}_n$ and $S\in {\cal P}_n$ by the isomorphism of the groups ${\cal P}_n$ and ${\cal S}_n$ described by the formulas (\ref{pi}), (\ref{sigma}), (\ref{P}) and (\ref{S}), then for each $t=1,2,\ldots ,n$ there are:
\begin{equation}\label{l21}
(t)\pi = t\oplus 1 \quad (\mbox{by definition})
\end{equation}
\begin{equation}\label{l22}
(t)\sigma = n\oplus 1\ominus t\quad (\mbox{by definition})
\end{equation}

We prove by induction, that for each whole positive number $l$ is true:
\begin{equation}\label{l23}
(t)\pi^l =t\oplus l
\end{equation}
When $l=1$, the proposition follows from (\ref{l21}). Let the equation $(t)\pi^l =t\oplus l$ be true. Then we get $(t)\pi^{l+1} =((t)\pi^l)\pi =(t\oplus l)\pi = t\oplus l\oplus 1$, and it follows that the equation (\ref{l23}) is true for every whole positive number $l$.

Using the equations (\ref{l21}), (\ref{l22}) и (\ref{l23}) we continuously get:
$$(t)\pi^l \sigma =((t)\pi^l )\sigma =(t\oplus l)\sigma =n\oplus 1\ominus (t\oplus l)=n\oplus 1\ominus t \ominus l$$
$$(t)\sigma\pi^{n-l} = ((t)\sigma )\pi^{n-l} = (n\oplus 1\ominus t)\pi^{n-l} =(n\oplus 1\ominus t)\oplus (n-l) = 2n\oplus 1\ominus t\ominus l = n\oplus 1\ominus t\ominus l$$

The last equation is true, because $2n\equiv n\equiv 0\; (\mbox{\rm mod}\; n)$. We see, that $(t)\pi^l \sigma =(t)\sigma\pi^{n-l}$ for every $t=1,2,\ldots ,n$ and therefore, $\pi^l \sigma =\sigma\pi^{n-l}$. We have in mind the isomorphism of the groups ${\cal S}_n$ and ${\cal P}_n$ it follows that the proposition in the lemma is true.

\hfill $\Box$

\begin{theor}\label{t1}
If $A\sim A^S$ and $B\sim A$ then $B\sim B^S$.
\end{theor}

Proof. Since $B\sim A$, then according to the lemma \ref{l1} there exist $k,l\in \{ 1,2,\ldots ,n\}$, such as $B=P^k AP^l$. Applying lemma \ref{l2} we get $BS=P^k AP^l S=P^k ASP^{n-l}$, and then follows, that $BS\sim AS$, i.e. according (\ref{8}) $B^S \sim A^S$. But $A^S \sim A\sim B$ and because of the transitiveness of the relation $\sim$ we get $B^S \sim B\Rightarrow B\sim B^S$.

\hfill $\Box$

Theorem \ref{t1} gives us a reason to make the following definition:
\begin{dfn}\label{samoogledalna}
Let  $A\in {\cal Q}_n$. We say that $A$  is a representative of {\bf self-mirrored image} (or {\bf mirror image to itself})  interweaving, if
$A\sim A^S .$
\end{dfn}

Let denote  the set $$\overline{{\cal M}_n} \subset \overline{{\cal Q}_n }$$ including all self-mirrored interweavings with repetition equal to $n$.

For random $A\in {\cal B}_n$ we define the operation
\begin{equation}\label{A^R}
A^R =
\left(
\begin{array}{cccc}
a_{1\, n} & a_{2\, n} & \cdots & a_{n\, n} \\
a_{1\, n-1} & a_{2\, n-1} & \cdots & a_{n\, n-1} \\
\vdots    & \vdots    &        & \vdots    \\
a_{1\, 1} & a_{2\, 1} & \cdots & a_{n\, 1} \\
\end{array}
\right)
=\left( A^S \right)^T =\left( A S\right)^T =S^T  A^T = S A^T
\end{equation}

In other words the matrix $A^R$ is received by rotating the matrix $A$ by ${\rm 90}^\circ$ counter clockwise.

Obviously,
$$\left( \left( \left( A^R \right)^R \right)^R \right)^R =A .$$

In the common case  $A^R \ne A$.

\begin{lemma}\label{l3}
If $P$ is a binary matrix, defined by the formula (\ref{P}), then
$$P^T =P^{n-1} $$
\end{lemma}

Proof. If $P=(p_{ij} )$ and $P^T =(p_{i j}' )$, then by definition $p_{i j}' = p_{ji}$ for each $i,j\in \{ 1,2,\ldots n\}$. Let $PP^T =Q=(q_{i,j} )\in {\cal P}_n$. Then for each $i=1,2,\ldots ,n$ there is $\displaystyle q_{ii} = \sum_{k=1}^n p_{ik} p_{ki}' =\sum_{k=1}^n p_{ik}^2 =(n-1)0 +1=1$ and it is the unique one in $i$th row of the matrix $Q=PP^T$. Therefore, $PP^T =E_n$, where $E_n$ is the identity matrix. We multiply from the left the two sides of the last equation with $P^{n-1}$ and have in mind, that $P^n =E_n$, then we get $P^{n-1} PP^T =P^{n-1} E_n$, and then finally we get, that $P^T =P^{n-1}$.

\hfill $\Box$

\begin{theor}\label{t2}
If $A\sim A^R$ and $B\sim A$ then $B\sim B^R$.
\end{theor}

Proof. $B\sim A$, hence according to the lemma \ref{l1} there exist natural numbers $k$ and $l$, as $B=P^k AP^l$. Then, when we apply lemma \ref{l2} and lemma \ref{l3} we get $B^R =SB^T =S(P^k A P^l )^T = S(P^T )^l A^T (P^T )^k =S(P^{n-1} )^l A^T (P^{n-1} )^k = SP^{nl-l} A^T P^{kn-k} = SP^{n-l}  A^T P^{n-k} =P^l S A^T P^{n-k}$. Therefore, $$B^R \sim A^R \sim A\sim B.$$
\hfill $\Box$

Theorem \ref{t2} gives us the right to give the following definition:
\begin{dfn}\label{rotat}
Let $A\in {\cal Q}_n$. If   $A \sim A^R$, then we say that $A$ is a representative of  {\bf rotation stable} interweaving.
\end{dfn}

Let denote   the set $$\overline{\cal R}_n \subset \overline{\cal Q}_n$$ of all rotation stable interweavings with repetition equal to $n$.

The rotation stable interweavings play important role in practice.This means, that if a fabric is weaved which weaving structure is coded with a matrix, representative of rotation stable interweaving, then this fabric will have the same operating characteristics (except of course the color) after a rotation by ${\rm 90}^\circ$.

\section{Quantity evaluation of the sets of all self-mirrored and all rotation stable interweavings with given repetition $n$.}
 In \cite{umb2010} is described a representation of the elements of ${\cal B}_n$ using ordered $n$-tuples of natural numbers $<k_1 ,k_2 ,\ldots ,k_n >$, where $0\le k_i \le 2^n -1$, $i=1,2,\ldots ,n$. The one to one corresponding is based on the definite representation of the natural numbers in binary number system, i.e. the number $k_i$ in binary number system (having eventually some zeroes in the beginning)  is the  $i$th row of the corresponding binary matrix. In \cite{umb2010} is proved that using this representation there are faster and saving memory algorithms. Having in mind this we create an algorithm, which finds just one representative of each equivalence class to the factor sets $\overline{{\cal Q}_n}$, $\overline{{\cal M}_n}$ and $\overline{{\cal R}_n}$. And the representative we receive is the minimal of the equivalence class with regard to the lexicographic order, this order is naturally brought in the set $\mathbb{N}^n$ of all ordered $n$-tuples of whole nonnegative numbers. Therefore, we get an algorithm to solve the combinatoric problem to find the number of the equivalence classes in the sets ${\cal Q}_n$, ${\cal M}_n$ and ${\cal R}_n$ relevant to the relation $\sim$ with given natural number $n$.

The matrices $P$ and $S$ given by the formulas (\ref{P}) and (\ref{S}) are coded using ordered $n$-tuples as it follows:
\begin{equation}\label{Ptuple}
P\; : \; <2^{n-2},2^{n-3} ,\ldots ,2^1 ,2^0 ,2^{n-1} >
\end{equation}
\begin{equation}\label{Stuple}
S\; : \; <2^0,2^1 , 2^2 ,\ldots  ,2^{n-2} ,2^{n-1} >
\end{equation}

In some programming languages (for example C, C++, Java \cite{umb2010,umb2009}) the number $x=2^k$ is calculated using ones the operation bitwise shift left $''<<''$ and the operator (statement) $x=1<<k;$.

In \cite{umb2010} is entered the operation \textit{logical multiplication} of two binary matrices, which we denote ''$*$''. Let $A=(a_{ij} )$ and $B=(b_{ij} )$ are matrices of ${\cal B}_n$. Then
\begin{equation}\label{logmult}
A*B =C=(c_{ij} )\in {\cal B}_n
\end{equation}
and by definition for each $i,j \in \{ 1,2\ldots ,n \} $
\begin{equation}\label{cij}
c_{ij} = \bigvee_{k=1}^n (a_{ik}  \; \& \; b_{kj} ),
\end{equation}
where we denote $\&$ and $\vee$ the operations conjunction and disjunction  in the boolean algebra ${\cal B}_n (\& ,\vee )$.

Analogously of the classical proof that the operation multiplication of matrices is associative (see for example. \cite{kurosh}) is proved, that the operation logical multiplication of binary matrices is associative. Therefore, ${\cal B}_n$ with the entered operation logical multiplication is monoid with identity - the identity matrix $E_n$. However, ${\cal P}_n$ is not trivial subgroup of this monoid.

 In \cite{umb2010,umb2009} is described an algorithm, needed $O(n^2 )$ operations and made the operation logical multiplication of two binary matrices, which are represented using ordered $n$-tuple. In the same time to make the product of two matrices according to the classical definition we need $O(n^3 )$ operations.

If the binary matrix $A$ is represented using the ordered $n$-tuple of numbers, then to check whether $A$ belongs to the set ${\cal Q}_n \subset {\cal B}_n$ we can use the following obvious proposition:

\begin{lemma}\label{l4}
 Let $A\in {\cal B}_n$ and let $A$ is represented by the ordered $n$-tuple $<k_1 ,k_2 ,\ldots ,k_n >$, where $0\le k_i \le 2^n -1$, $i=1,2,\ldots ,n$ and let denote $|$ and $\&$ the operations bitwise ''or'' and bitwise ''and'' (to get detailed definitions see for example \cite{davis}, \cite{Kernigan} or \cite{umb2009}). Then:

(i) The number $k_i$, $i=1,2,\ldots ,n$ represents row of zeroes if and only if $k_i =0$;

(ii) The number $k_i$, $i=1,2,\ldots ,n$ represents row of ones if and only if $k_i =2^n -1$;

(iii) $j$th column of $A$ is made of zeroes if and only if $$(k_1 \, | \, k_2 \, | \cdots |\, k_n )\, \& \, 2^j =0;$$

(iv) $j$th column of $A$ is made of ones if and only if $$(k_1 \, \& \, k_2 \, \& \cdots \& \, k_n )\, \& \, 2^j \ne 0.$$
\hfill $\Box$
\end{lemma}

The algorithm, which is below-described is based on the following propositions:
\begin{lemma}\label{l5}
If $A\in {\cal B}_n $, $B\in {\cal P}_n$, then
$$A*B=AB$$ and
$$B*A=BA$$
\end{lemma}

Proof. Let $A=(a_{ij} )$, $B=(b_{ij} )$, $U=A*B=(u_{ij} )$ and $V=AB=(v_{ij} )$, $i,j=1,2,\ldots ,n$. Let the unique one in the $j$th column of $B\in {\cal P}_n$ is on the $s$th place, i.e. $b_{sj} =1$ and $b_{kj} =0$ when $k\ne s$. Then by definition
$$u_{ij} =\bigvee_{k=1}^n (a_{ik}  \; \& \; b_{kj} )=
\left\{
\begin{array}{ccc}
1 & \mbox{\rm for}  & a_{is} =1\\
0 & \mbox{\rm for}  & a_{is} =0
\end{array}
\right.$$ and
$$v_{ij} =\sum_{k=1}^n (a_{ik}  b_{kj} )=
\left\{
\begin{array}{ccc}
1 & \mbox{\rm for}  & a_{is} =1\\
0 & \mbox{\rm for}  & a_{is} =0
\end{array}
\right.$$

Therefore, $u_{ij} = v_{ij}$ for each $i,j\in \{ 1,2,\ldots , n\}$.

Analogously can be proved, that $B*A=BA$.

\hfill $\Box$

\begin{lemma}\label{l6}
 Let $A\in {\cal B}_n$ is represented using ordered $n$-tuple $<k_1 ,k_2 ,\ldots ,k_n >$ and let $A$ is a minimal element of the equivalence class corresponding to the lexicographic order in $\mathbb{N}^n$. Then $k_1 \le k_t$ for each $t=2,3,\ldots ,n$.
\end{lemma}

Proof. We presume, that there exists $t\in \{ 2,3,\ldots ,n\}$, such as $k_t <k_1$. Then if we move the first row on the last place $t-1$ times, we get a matrix $A' \in{\cal B}_n$, such as $A' \sim A$ and $A'$ is represented using the $n$-tuple $<k_t ,k_{t+1} ,\ldots ,k_n ,k_1 ,\ldots ,k_{t-1} >$. Then obvious $A' <A$ according to the lexicographic order in $\mathbb{N}^n$, which runs counter to the minimum of $A$ in the equivalence class $\overline{A}$.

\hfill $\Box$

We can create the following generalized algorithm to obtain just one representative of each equivalence class in the factor sets $\overline{\cal Q}_n$, $\overline{\cal M}_n$ and  $\overline{\cal R}_n$

\begin{algr}\label{alg1}.
\begin{enumerate}
\item\label{it1} Generating all ordered $n$-tuples of natural numbers $<k_1 ,k_2 ,\ldots ,k_n >$ such as $1\le k_i \le 2^n -2 $ and $k_1 \le k_i$ as $i=2,3,\ldots ,n$;
\item\label{it2} Check if the elements obtained in \ref{it1} belong to the set ${\cal Q}_n$ according to the lemma \ref{l4} (iii)  (iv) (cases (i) and (ii) we reject when we generate the elements in point \ref{it1} according to the lemma \ref{l6});
\item\label{it3} Check whether the element, obtained in point \ref{it2} is minimal in the equivalence class. According to the lemmas \ref{l1} and \ref{l5} $A$ is minimal in $\overline{A}$ if and only if $A\le P^k *A*P^l$ for each $k,l\in \{ 0,1,\ldots ,n-1\}$, where $P$ is the matrix represented by $n$-tuple (\ref{Ptuple});
\item\label{it4} Check whether the elements obtained in point \ref{it3} belong to the set ${\cal M}_n$ according to definition \ref{samoogledalna} and applying lemmas \ref{l1} and \ref{l5};
\item\label{it5} Check whether the elements obtained in point \ref{it3} belong to the set ${\cal R}_n$ according to definition  \ref{rotat} and applying lemmas \ref{l1} and \ref{l5}.
\end{enumerate}
\hfill $\Box$
\end{algr}

The work results of the algorithm \ref{alg1} taking some values of $n$ are generalized in the following table:

\begin{center}
\begin{tabular}{||c||c|c|c|c|c|c||}
\hline
$n$                                      &   2  &   3   &   4       &      5    \\
\hline
$\displaystyle |\overline{\cal Q}_n |$   &   1  &   14  &  1 446    & 705 366   \\
\hline
$\displaystyle |\overline{\cal M}_n |$   &   1  &   2   &    142    &   1 302   \\
\hline
$\displaystyle |\overline{\cal R}_n |$   &   1  &   2   &     18    &      74   \\
\hline
\end{tabular}
\end{center}
When $n\ge 6$ are got too large values (see (\ref{bn})) and to avoid ''overloading'' it is necessary to be used some special programming techniques which is not the task in this work.

$
\begin{array}{llllll}
\mbox{Krasimir Yankov Yordzhev}   & & & & & \mbox{Hristina Aleksandrova Kostadinova}\\
\mbox{South-West University ''N. Rilsky''} & & & & & \mbox{South-West University ''N. Rilsky''}\\
\mbox{2700 Blagoevgrad, Bulgaria} & & & & & \mbox{2700 Blagoevgrad, Bulgaria}\\
\mbox{Email: yordzhev@swu.bg, iordjev@yahoo.com} & & & & & \mbox{Email: hkostadinova@gmail.com}
\end{array}
$


\begin {thebibliography}{99}
\bibitem{davis} \textsc{S. R. Davis} C++ for dummies.
\emph{IDG Books Worldwide,} 2000.
\bibitem{Kernigan} \textsc{B. W. Kernigan, D. M Ritchie}  The C programming Language.
\emph{ AT$\&$T Bell Laboratories,}  1998.
\bibitem{umb2010} \textsc{H. Kostadinova, K. Yordzhev} A Representation of Binary Matrices
\emph{Mathematics and education in mathematics},  39 (2010), 198--206.
\bibitem{curtis}
\textsc{C. W. Curtis, I. Rainer} Representation Theory of Finite Groups and Associative Algebras.
\emph{ Wiley-Interscience}, 1962.
\bibitem{hall}
\textsc{M. Hall} The Theory of Groups.
\emph{ New York}, 1959.
\bibitem{huppert}
\textsc{B. Huppert} Enlishe Gruppen.
\emph{Springer}, 1967.

\bibitem{bansko} \textsc{K. Yordzhev} On an equivalence relation in the set of the permutation matrices.
\emph{Blagoevgrad, Bulgaria, SWU, Discrete Mathematics and Applications},  (2004), 77--87.

\bibitem{umb2009} \textsc{K. Yordzhev} An example for the use of bitwise operations
in programming.
\emph{Mathematics and education in mathematics},  38 (2009), 196--202.

\bibitem{borzunob} \textsc{Г. И. Борзунов } Шерстяная промишленост-обзорная информащия.
\emph{Москва, ЦНИИ ИТЭИЛП,}  3, (1983).

\bibitem{yordzhev} \textsc{К. Я. Йорджев, И. В. Статулов}  Математическо моделиране и количествена оценка на първичните тъкачни сплитки.
\emph{ Текстил и облекло, }  10, (1999), 18--20.

\bibitem{kurosh} \textsc{А. Г. Курош}  Курс высшей алгебры.
\emph{Москва, Наука,}  1975.

\bibitem{melnikov}
\textsc{О. В. Мельников, В. Н. Ремесленников, В. А. Романков, Л. А. Скорняков, И. П. Шестяков} Общая алгебра.
\emph{ Москва, Наука}, 1990.


\bibitem{tarakanov} \textsc{В. Е. Тараканов} Комбинаторные задачи и (0,1)-матрицы.
\emph{Москва, Наука,}  1985.


\end{thebibliography}
\end{document}